\documentclass[preprint,3p,times]{elsarticle}

\usepackage[english]{babel}
\usepackage{amsthm}
\usepackage{amsmath}
\usepackage{amsfonts}
\usepackage{amssymb}
\usepackage{makeidx}
\usepackage{graphicx}
\usepackage{hyperref}
\usepackage{xcolor}
\usepackage[justification=centering]{caption}
\usepackage{subcaption}
\usepackage{algorithm}
\usepackage[noend]{algpseudocode}
\usepackage{cleveref}
\crefname{equation}{Eq.}{Eqs.}
\usepackage{balance}
\usepackage{mathtools}

\DeclareMathOperator*{\Argmax}{Arg\,max}

\newcommand{\markupdraft}[2]{
	\ifthenelse{\equal{#1}{display}}{#2}{}
	\ifthenelse{\equal{#1}{color}}{\color{#2}}{}
}

\newcommand{\newcolored}[3][]{{\markupdraft{color}{#2}#3}
	\ifthenelse{\equal{#1}{}}{}{\markupdraft{display}{{\color{yellow!70!black}[#1]}}}} 

\newcommand{\del}[2][]{{\markupdraft{display}{{\color{yellow!99!black}[removed: "#2"[#1]]}}}} 
\newcommand{\dell}[2][]{{\markupdraft{display}{{\color{gray}[removed: "#2"#1]}}}} 
\newcommand{\new}[2][]{\newcolored[#1]{blue!75!black}{#2}}

\newcommand{\indraftonly}[1]{{#1}}  
\renewcommand{\indraftonly}[1]{}\renewcommand{\markupdraft}[2]{}  
\renewcommand{\del}[2][]{}

\newcommand{\blank}[1]{}

\newcommand{\vertiii}[1]{{\left\vert\kern-0.25ex\left\vert\kern-0.25ex\left\vert #1 
		\right\vert\kern-0.25ex\right\vert\kern-0.25ex\right\vert_2}}

\newcommand*{\1}{\text{\usefont{U}{bbold}{m}{n}1}}

\newcommand{\bR}{\mathbb{R}}

\newcommand{\bN}{\mathbb{N}}

\newcommand{\cO}{\mathcal{O}}

\newcommand{\diag}{\mathrm{diag}}

\newtheorem{proposition}{Proposition}
\newtheorem{theorem}[proposition]{Theorem}

\newtheorem{lemma}[proposition]{Lemma}

\usepackage{xparse}

\NewDocumentCommand{\yt}{ O{t} }{Y_{#1}}
\NewDocumentCommand{\yit}{ O{i} O{t} }{\left[Y_{#2}\right]_{#1}}

\newcommand{\maxtwo}[2]{\max\{#1,#2\}}
\newcommand{\mintwo}[2]{\min\{#1,#2\}}

\newcommand{\fracp}[2]{#1/#2}

\usepackage{enumitem} 
\setlist[enumerate]{wide = 0pt, leftmargin=*}

\newcommand{\shorter}[2][]{#1}

\usepackage{caption, subcaption}

\usepackage{enumitem} 
\setlist[enumerate]{wide = 0pt, leftmargin=*}

\journal{Applied Mathematics Letters}

\renewcommand{\ge}{\geqslant}

\begin{document}
	
	\begin{frontmatter}
		
		
		
		\title{Asymptotic estimations of a perturbed symmetric eigenproblem 
		}

		\cortext[cor1]{Corresponding author}
		
		\author[cmap]{Armand Gissler\corref{cor1}}
		\ead{firstname.lastname@polytechnique.edu}
		
		\affiliation[cmap]{organization={Inria, CMAP, CNRS, École polytechnique, Institut Polytechnique de Paris},
			city={91120 Palaiseau},
			country={France}}

		\author[cmap]{Anne Auger}
		\ead{firstname.lastname@inria.fr}
		\author[cmap]{Nikolaus Hansen}
		\ead{firstname.lastname@inria.fr}
		
		
		\begin{abstract}
			We study ill-conditioned positive definite matrices that are disturbed by the sum of $m$ rank-one matrices of a specific form.
			We provide estimates for the eigenvalues and eigenvectors.
			When the condition number of the initial matrix tends to infinity,
			we bound the values of the coordinates of the eigenvectors of the perturbed matrix.
			Equivalently, in the coordinate system where the initial matrix is diagonal, we bound the rate of convergence of coordinates that tend to zero.
		\end{abstract}
		
		
		
		
	\end{frontmatter}
	
	
	\section{Introduction}

Given a $d\times d$ symmetric matrix 
with known eigenvectors and eigenvalues denoted $B$, and a rank-one matrix $vv^T$ where $v \in\bR^d$, eigenvalues and eigenvectors of matrices of the form 
	\begin{equation}
		A = B + vv^T
		\label{eq:sym-plus-rank1}
		\tag{$P_1$}
	\end{equation}
	have been widely studied, notably in the context of perturbation theory.
	For instance, the eigenvalues of \eqref{eq:sym-plus-rank1}
	can be estimated and a 
	formula for the eigenvectors is known~\cite{golub1973modified,bunch1978rank,ipsen2009refined}.
		Specifically, if $B$ is a diagonal matrix $\diag(\lambda_1,\lambda_2,\dots,\lambda_d)$ with distinct eigenvalues and $v$ has only nonzero entries, then the component $j$ of the unit eigenvector associated to eigenvalue  $\nu_i$ of the updated matrix $A$ satisfies the so-called Bunch-Nielsen-Sorensen formula
		\begin{equation}
			 C_i \times \frac{[v]_j}{\lambda_j-\nu_i} \quad \text{for~}i,j=1,\dots,d
			\label{eq:BNS-formula}
		\end{equation}
		where $C_i$ is a nonzero normalization constant
		and $[\,.\,]_j$ denotes the $j$-th coordinate, a notation we will continue to use in the sequel.
Several results have been established for additive perturbations of rank $1$ \cite{ding2007eigenvalues,benasseni2019inequalities} and of higher rank \cite{thompson1976behavior,mathias1997spectral}.  
	Symmetric and nonsymmetric perturbation eigenvalue problems have been studied \cite{bhatia2007perturbation} 
	as well as perturbation \new{results for invariant subspaces} \cite{karowperturbation}.
	In this paper, we provide relative perturbation bounds for the eigenvectors of positive definite matrices.
	In contrast to previous relative perturbation results for eigenvalues~\cite{ipsen1998relative} and invariant subspaces~\cite{truhar1999relative,ipsen2000overview},
	the bounds in our result depend 
	on the eigenvalues of
	the initial matrix $B$ rather than the norm of the perturbation, see \Cref{eq:rank-mu-eigenvectors-symmetric} below.

	
	Specifically, we consider the perturbation with a sum of $m$ rank-one matrices of the form
	\begin{equation}
		A^{(m)} = B + \sqrt{B}\sum_{i=1}^m [v^{(i)}] [v^{(i)}]^T \sqrt{B}
		\label{eq:diagonal-plus-rankm}
		\tag{$P_m$}
	\end{equation}
	with $B=P\diag(\lambda_1,\dots,\lambda_d)P^T$\shorter{ a positive symmetric matrix} where $P$ is an orthogonal matrix, $\lambda_1\geqslant\dots\geqslant\lambda_d>0$ are the\shorter{ (decreasingly ordered)} eigenvalues of $B$, and $v^{(1)},\dots,v^{(m)}\shorter[\in]{$ are vectors of $}\bR^d$.
	The square root $\sqrt{B}:=P\diag(\sqrt{\lambda_1},\dots,\sqrt{\lambda_d})P^T$ is defined as the unique symmetric positive definite matrix such that $\sqrt{B}\times\sqrt{B}=B$, see e.g.\ \cite[Theorem 7.2.6]{horn2013matrix}.
	Matrices of the form \eqref{eq:diagonal-plus-rankm} are used in various applications in different domains. For instance, low rank updates of covariance matrices are used in stochastic optimization~\shorter[\cite{kjellstrom1981stochastic,hansen2003reducing}]{\cite{kjellstrom1981stochastic}}, system identification~\cite[p.\,369]{ljung1999system},\shorter[ and]{} adaptive Markov Chain Monte Carlo methods~\cite{haario2001adaptive}\shorter{, and rank-$\mu$ updates, with $\mu\in\bN^*$, for covariance matrices in the stochastic optimization algorithm CMA-ES \cite{hansen2003reducing}}.
	Our motivation is to study the eigenvectors of $A$ in \eqref{eq:diagonal-plus-rankm}, denoted as $e_i^{(m)}$ in the sequel,
	\emph{when the matrix $B$ is highly ill-conditioned}\shorter{ (i.e., when $\max_{i} \lambda_i/\min_{i}\lambda_i$ is large)}.
	When $d=2$ and $m=1$ we can compute the eigenvectors of $A^{(1)}$ explicitly. As an example, consider $B = \diag(\lambda_1,1)$ where $\lambda_1 >1$ and $v^{(1)}=[1,1]^T$.
	Then, the unit eigenvector associated to the largest eigenvalue\shorter{ the vector $e_1^{(1)}$ where} of $A^{(1)}$
	obeys $\sqrt{1+s^2} \times e_1^{(1)} = [1 , s]^T$
	with $s=\lambda_1^{1/2}\times(1-\lambda_1^{-1}-\sqrt{1-\lambda_1^{-1}+\lambda_1^{-2}})=-\lambda_1^{-1/2}/2+\cO(\lambda_1^{-3/2})$
	and hence $[e_1^{(1)}]_2=\lambda_1^{-1/2} + o(\lambda_1^{-1/2})$ when  $\lambda_1\to\infty$.
	Hence, the (second) coordinate of the (first) unit eigenvector of $A^{(1)}$ vanishes like $1/\sqrt{\lambda_1}$ when $\lambda_1\to\infty$. In this paper, we generalize this result to the case where $d \geqslant 2$ and $m \geqslant 1$,
	as summarized in the following theorem which directly follows from \Cref{th:eigenvectors-rankm} below.

	\begin{theorem}
		If $e_i^{(m)}$ is a unit eigenvector corresponding to the $i$-th largest eigenvalue (counted with multiplicity)
		of $A^{(m)}$ in \eqref{eq:diagonal-plus-rankm}
		and $e_j^{(0)}$ is a unit eigenvector corresponding to the $j$-th largest eigenvalue $\lambda_j$ of $B$, then
		\begin{equation}
			\left|\left\langle e_i^{(m)}, e_j^{(0)} \right\rangle\right| \leqslant C_m\times\sqrt{\frac{\mintwo{\lambda_i}{\lambda_j}}{\maxtwo{\lambda_i}{\lambda_j}}}
			\label{eq:rank-mu-eigenvectors-symmetric}
		\end{equation}
		where $C_m>0$ is a constant which depends polynomially on $d$ and $\max_{k=1,\dots,m} \| v^{(k)}\|$.
		\label{th:eigenvectors-rankm-symmetric}
	\end{theorem}
When $B$ is diagonal, $e_j^{(0)}$ is the $j$-th canonical unit vector. Hence $|\langle e_i^{(m)}, e_j^{(0)} \rangle | = [e^{(m)}_i]_j$ and the theorem implies in particular that the $j$-th coordinate of $e^{(m)}_i$ converges to zero at least as fast as
$\sqrt{\mintwo{\lambda_i}{\lambda_j}/\maxtwo{\lambda_i}{\lambda_j}}$ when the latter tends to $0$
(which is tight in the above example when $d=2$ and $m=1$),
thereby limiting the change of the angle between these eigenvectors. 
Considerations on the angle between eigenspaces have been made previously~\cite{davis1970rotation}, however matrices on the form of \eqref{eq:diagonal-plus-rankm} have not been studied in this context.
In the remainder, 
we always choose w.l.o.g.\ 
the coordinate system where the matrix $B$ of \eqref{eq:diagonal-plus-rankm} is diagonal and has decreasingly ordered diagonal values.

	This inequality is crucial to study 
	the stability of a Markov chain underlying the CMA-ES algorithm~\cite{hansen2001completely,hansen2003reducing}.
	Proofs of linear convergence for Evolutionary Strategies (ES) rely on a drift condition \cite[Theorem 17.0.1]{meyn2012markov} to prove the ergodicity of an underlying Markov chain, see e.g.\ \cite{auger2016linear,toure2023global}.
	To apply this approach to CMA-ES, a potential function is defined on the state-space of this Markov chain and its expected decrease is proven outside a compact set.
	The state space includes a covariance matrix, 
    updated as
	\begin{equation}\label{eq:CMA-rankmu}
		C_{t+1} = (1-c) C_t + c \sqrt{C_t} \sum_{i=1}^m w_i U_{i}U_{i}^T \sqrt{C_t}\enspace,
	\end{equation}
	where $c\in[0,1]$, $w_1,\dots,w_m$ are positive weights that sum to $1$, and the vectors $U_{i}$, $i=1,\dots,m$, are Gaussian vectors ranked according to a fitness function \cite[Eq.\;(11)]{hansen2003reducing}.
	Hence, \eqref{eq:diagonal-plus-rankm} encompasses the update of this covariance matrix.
	  \Cref{eq:rank-mu-eigenvectors-symmetric}
	  is needed to bound
	  the expected condition number of the updated covariance matrix, since it controls 
	  the influence of small eigenvalues on the growth of the largest eigenvalues.
  
	%
	%

	This paper is organized as follows. In \Cref{sec:eigenvalues}, we study the eigenvalues of \eqref{eq:diagonal-plus-rankm}. 
	 In \Cref{sec:eigenvectors}, we provide bounds for the coordinates of the eigenvectors using \Cref{eq:BNS-formula}, and provide an empirical result suggesting that these bounds are tight.

	\section{Bounds on the eigenvalues of \eqref{eq:diagonal-plus-rankm}}\label{sec:eigenvalues}
	
	The Bunch-Nielsen-Sorensen formula \eqref{eq:BNS-formula}
	which we will use in \Cref{sec:eigenvectors}
	requires the eigenvalues of the \emph{updated} matrix.
	Thus,
	we first derive bounds
	on the (decreasingly ordered) eigenvalues
	\begin{equation}
		\label{eq:minmax-principle}
		\lambda_i (A^{(m)}) = \max_{V\subset \bR^d, \dim V=i} \,\min_{v\in V, v\neq0} \frac{v^TA^{(m)}v}{v^Tv} 
		= \min_{V\subset \bR^d, \dim V=d-i+1}
		\,\max_{v\in V, v\neq0} \frac{v^TA^{(m)}v}{v^Tv}
		\quad\text{for~} i=1,\dots,d
	\end{equation}	
	where the equalities ensue from the min-max principle
	and from Gersgorin's circle theorem
	\cite[Theorems~4.2.6 and 6.1.1]{horn2013matrix}.
	
	Theorem 2.7 in \cite{ipsen1998relative} and Theorem 2.1 in \cite[p.~175]{stewart2001matrix} provide an estimation for the eigenvalues of \eqref{eq:diagonal-plus-rankm}:
	\begin{equation}\label{l:eigenvalues-rankm}
		\lambda_i \leqslant \nu_i \leqslant \lambda_i \times \left(1+md\times\max_{k=1,\dots,m}\|v^{(k)}\|_\infty^2\right)
		\quad\text{for~}i\in\{1,\dots,d\}\enspace.
	\end{equation}
	The next lemma provides a slightly tighter
	upper bound on these eigenvalues after a single rank-one pertubation
	($m=1$) and is used in \Cref{p:eigenvectors-of-rank-one-update-distincts-eigenvalues}.
%
	
	\begin{lemma}
		\label{l:eigenvalues-rank1}
		Let $D=\diag(\lambda_1,\dots,\lambda_d)$ be a diagonal matrix with $\lambda_1>\dots>\lambda_d>0$. Let $v\in\bR_{\neq 0}^d$ be a vector with only nonzero entries. Let $A=D+\sqrt{D}vv^T\sqrt{D}$ and
		$\nu_1\geqslant\nu_2\geqslant\dots\geqslant\nu_d$
		denote
		the eigenvalues of $A$. Then,
%
			\begin{equation}
				\nu_i \leqslant \lambda_i \times \left(1+(d-i+1)\| v \|_\infty |[v]_{j_i}|\right)
				\quad\text{for all~}i\in\{1,\dots,d\}
				\label{eq:majoration-eigvalue-rankone}
			\end{equation}
			where
			$
				j_i\in\Argmax_{j=i,\dots,d}\left\{|[v]_j|
				\;\big|\; 
				\lambda_j \geqslant \lambda_i\times \left(1-\sqrt{\frac{\lambda_j}{\lambda_i}}(d-i+1)\| v \|_\infty |[v]_j|\right) \right\}
				\label{eq:def-j_i}
				$.
	\end{lemma}
	\begin{proof}
		%
		Fix $i\in\{1,\dots,d\}$ and remark that, by \Cref{eq:minmax-principle}, we have
		$
		\nu_i \leqslant  \max_{v \in \bar{V}_i,\|v\|=1} {v^T Av} =  \lambda_1\left( [A]_{i:d,i:d} \right) ,
		$
		where $\bar{V}_i = {\rm Vect}(e_{i},\ldots,e_{d}) $ with $e_i$ being the $i^{\rm th}$ vector of the standard basis of $\bR^d$, and with $[A]_{i:d,i:d}$ denoting the submatrix of $A$ from rows and columns with indices between $i$ and $d$ included.
		But, by \cite[Theorem 6.1.1]{horn2013matrix}, we also have that
		$$
		\lambda_1\left( [A]_{i:d,i:d} \right) \leqslant \max_{j=i,\dots,d} \left(  \sum_{k=i}^d |[A]_{j,k}|  \right)  \eqqcolon \max_{j\geqslant i} B_j.
		$$
		Since $A=D+\sqrt{D}vv^T\sqrt{D}$, then $|[A]_{j,k}|\leqslant \sqrt{\lambda_j\lambda_k}(\1\{j=k\}+\|v\|_\infty|[v]_j|)$. If $j\geqslant i$ is such that $\lambda_j\geqslant \lambda_i\times (1-\sqrt{\fracp{\lambda_j}{\lambda_i}}(d-i+1)\|v\|_\infty |[v]_j|)$, then by definition of $j_i$ we have then $|[v]_j|\leqslant |[v]_{j_i}|$, yielding to $B_j\leqslant \lambda_i\times (1+ (d-i+1)\|v\|_\infty|[v]_{j_i}|)$. Any other $j\geqslant i$ satisfies
		$
			\lambda_j < \lambda_i - \sqrt{\lambda_j\lambda_i}(d-i+1)\|v\|_\infty |[v]_j| ,
		$
		hence by sum $B_j\leqslant \lambda_i$. All in all, $\max_{j\geqslant i} B_j \leqslant \lambda_i\times (1+ (d-i+1)\|v\|_\infty|[v]_{j_i}|)$, proving \Cref{eq:majoration-eigvalue-rankone}. \dell{As $|[v]_{j_i}|\leqslant\|v\|_\infty$, we get \Cref{eq:majoration-eigvalue-rankonev2dell}.}
	\end{proof}
%

	\section{Estimating the eigenvectors of \eqref{eq:diagonal-plus-rankm}} \label{sec:eigenvectors}
	
We use the bounds from \Cref{l:eigenvalues-rank1,l:eigenvalues-rankm}
to estimate the eigenvectors of \eqref{eq:diagonal-plus-rankm}
by applying \Cref{eq:BNS-formula},
for $m=1$ in the next section and $m\ge1$ in \Cref{sec-rankm}.
	
	\subsection{Rank-one perturbation} \label{sec-rank1}
	
In \Cref{p:eigenvectors-of-rank-one-update-distincts-eigenvalues}, we obtain bounds on the coordinates of the eigenvectors of \eqref{eq:diagonal-plus-rankm} ($B$ is assumed to be diagonal)
when $m=1$, which comes as a consequence of \Cref{eq:BNS-formula}.

	\begin{proposition}
		\label{p:eigenvectors-of-rank-one-update-distincts-eigenvalues}
		Let $D=\diag(\lambda_1,\dots,\lambda_d)$ be a diagonal matrix with $\lambda_1\geqslant\dots\geqslant\lambda_d>0$.
		Let $v\in\bR^d$ and $V := \maxtwo{d^{-1/2}}{\|v\|_\infty}$. Consider the matrix $A=D+\sqrt{D}vv^T\sqrt{D}$ and $\nu_1\geqslant\dots\geqslant\nu_d$ its eigenvalues
		and 
		$(e_1^{(1)},\dots,e_d^{(1)})$
		a corresponding
		orthonormal basis of eigenvectors.
        Then,
		
		%
		\begin{equation}
			\left|[e_i^{(1)}]_j\right|\leqslant 5d^2V^4\sqrt{\frac{\mintwo{\lambda_i}{\lambda_j}}{\maxtwo{\lambda_i}{\lambda_j}}}
			\quad\text{for all~} i,j\in\{1,\dots,d\}
			\label{eq:rankone-eigenvectors}
		\end{equation}
	\end{proposition}
	
	\begin{proof}
		We prove first that, if $\maxtwo{\lambda_i}{\lambda_j}>(1+dV^2)\times\mintwo{\lambda_i}{\lambda_j}$, then
		\begin{equation}
			\label{eq:eigenvectors-rankone-lambdai>lambdaj-distinct}
 \left|[e_i^{(1)}]_j\right| \leqslant (d-i+1)V^2 \times  \frac{\inf_{\rho\in(0,1)}\psi(\rho,(d-i+1)V^2) }{1-(1+(d-i+1)V^2)\frac{\mintwo{\lambda_i}{\lambda_j}}{\maxtwo{\lambda_i}{\lambda_j}}}\sqrt{\frac{\mintwo{\lambda_i}{\lambda_j}}{\maxtwo{\lambda_i}{\lambda_j}}} 
		\end{equation}
		with $\psi(\rho,W)=\maxtwo{ 2(1-\rho)^{-1/2}}{2\rho^{-1} W } $, from which we deduce \Cref{eq:rankone-eigenvectors}.
		
		First suppose that the eigenvalues $\lambda_i$ of $D$ are distinct, and that all entries of $v$ are nonzero.
		Then, by \Cref{eq:BNS-formula}, we have for $i,j\in\{1,\dots,d\}$ that
		$
		[e_i^{(1)}]_j = C_i \frac{[\sqrt{D}v]_j}{\lambda_j-\nu_i}=C_i \frac{\sqrt{\lambda_i}[v]_j}{\lambda_j-\nu_i},
		$
		where $C_i\in\bR$ is chosen such that $\|e_i^{(1)}\|=1$, hence
		\begin{align}
			|C_i| & = \left\| \left( \frac{\sqrt{\lambda_j}[v]_j}{\lambda_j-\nu_i} \right)_{j=1,\dots,d} \right\|^{-1}= \left( \sum_{j=1}^d \left| \frac{\sqrt{\lambda_j}[v]_j}{\lambda_j-\nu_i} \right|^2 \right)^{-1/2}\leqslant \min_{1\leqslant j\leqslant d} \frac{|\lambda_j-\nu_i|}{\sqrt{\lambda_j}|[v]_j|} .
			\label{eq:maj-Ci-step1}
		\end{align}
		Combining \cite[Theorem 2.1, p.~175]{stewart2001matrix} with \Cref{eq:majoration-eigvalue-rankone}, we have
		$
			\lambda_i< \nu_i \leqslant \lambda_i\times (1+(d-i+1)V|[v]_{j_i}|),
		$
		where $j_i$ is defined in \Cref{eq:def-j_i}. By definition of $j_i$ we have
		$
			0 \leqslant \lambda_i -\lambda_{j_i} \leqslant \lambda_i (d-i+1) V |[v]_{j_i}| .
		$
		By sum, we obtain
		$
			0 <  \nu_i - \lambda_{j_i} \leqslant 2\lambda_i (d-i+1) V |[v]_{j_i}| .
		$
		We apply this to \Cref{eq:maj-Ci-step1} to get
		\begin{equation}
			|C_i| \leqslant \frac{\nu_i-\lambda_{j_i}}{\sqrt{\lambda_{j_i}}|[v]_{j_i}|} \leqslant \frac{2 \lambda_i (d-i+1)V|[v]_{j_i}|}{\sqrt{\lambda_{j_i}}|[v]_{j_i}|} = \frac{\lambda_i}{\sqrt{\lambda_{j_i}}} 2(d-i+1)V
			\enspace.
			\label{eq:maj-Ci-step2}
		\end{equation}
		Let $\rho\in(0,1)$. If $(d-i+1)V^2 \sqrt{\lambda_{j_i}} \leqslant \rho \sqrt{\lambda_i}$, then by definition of $j_i$, $\lambda_{j_i}\geqslant \lambda_i \times (1-\rho)$, and by \Cref{eq:maj-Ci-step2}, $|C_i|\leqslant (1-\rho)^{-1/2} \times 2(d-i+1)V\sqrt{\lambda_i}$. Otherwise, $|C_i|\leqslant \rho^{-1} \times 2(d-i+1)^2 V^3 \sqrt{\lambda_i}$. All in all, for $\rho\in(0,1)$,
		\begin{equation}
			|C_i| \leqslant (d-i+1)V\sqrt{\lambda_i} \times \max\left\{ 2(1-\rho)^{-1/2},2\rho^{-1} (d-i+1)V^2  \right\} \eqqcolon C_\rho \sqrt{\lambda_i}
			\enspace.
			\label{eq:maj-Ci-step3}
		\end{equation}
		Then, 
		$
		| [e_i^{(1)}]_j | = |C_i | \sqrt{\lambda_j}{|[v]_j|}/{|\lambda_j-\nu_i|} \leqslant  {\sqrt{\lambda_i\lambda_j}}/{|\lambda_j-\nu_i|} \times \inf_{\rho\in(0,1)}C_\rho.
		$
		By \Cref{eq:majoration-eigvalue-rankone}, when $\lambda_j<\lambda_i$,
		$
		|[e_i^{(1)}]_j| \leqslant  \min_{\rho\in[0,1]}C_\rho \times (1-\lambda_j/\lambda_i)^{-1}\sqrt{{\lambda_j}/{\lambda_i}} .
		$
		By \Cref{l:eigenvalues-rankm}, when $\lambda_j>(1+dV^2)\lambda_i$,
		$
		|[e_i^{(1)}]_j| \leqslant \min_{\rho\in[0,1]}C_\rho \times (1-(1+dV^2)\lambda_i/\lambda_j)^{-1}\sqrt{{\lambda_i}/{\lambda_j}} .
		$
		
		If the eigenvalues of $D$ are not distinct or not all entries of $v$ are nonzero,
		we consider a sequence of diagonal matrices 
		$
		\{D_k = \diag(\lambda_1^k,\dots,\lambda_d^k) \}_{k\in\bN}
		$
		such that the diagonal elements $\lambda_1^k>\lambda_2^k>\dots>\lambda_d^k>0$ are distinct and $D_k\to D$ when $k\to\infty$,
		and a sequence of vectors $\{v_k\in\bR^d_{\neq 0}\}_{k\in\bN}$ with only nonzero entries where $\|v_k\|_\infty\leqslant N$ and $v_k\to v$ when $k\to\infty$. Denote then $A_k = D_k +\sqrt{D_k}v_kv_k^T\sqrt{D_k}$ and $\nu_1^k\geqslant\dots\geqslant\nu_d^k$ its eigenvalues.
		Note that $A_k\to A$ when $k\to\infty$, so by continuity of the eigenvalues, $\nu_i^k\to \nu_i$ when $k\to\infty$. Furthermore, we just proved that if $e_1^k,\dots,e_d^k$ are unit eigenvectors of $A_k$ corresponding respectively to the eigenvalues $\nu_1^k,\dots,\nu_d^k$, then $e_1^k,\dots,e_d^k$ and $\lambda_1^k,\dots,\lambda_d^k$ satisfy \Cref{eq:eigenvectors-rankone-lambdai>lambdaj-distinct}. 
		Moreover, the vectors $e_i^k$ all belong to the unit sphere of $\bR^d$, so up to considering a subsequence of $\{A_k\}_{k\in\bN}$, we can assume w.l.o.g.\ that each $e_i^k$ tends to a vector $e_i^{(1)}\in\bR^d$ when $k\to\infty$. As $(e_1^k,\dots,e_d^k)$ is an orthonormal system of $\bR^d$, so is its limit $(e_1^{(1)},\dots,e_d^{(1)})$ and $e_i^{(1)}$ is an eigenvector of $A$ corresponding to the eigenvalue $\nu_i$. Therefore, \Cref{eq:eigenvectors-rankone-lambdai>lambdaj-distinct} holds by taking the limit $k\to\infty$ in the equation 
		satisfied by $e_1^k,\dots,e_d^k$ and $\lambda_1^k,\dots,\lambda_d^k$.
		
		To obtain \Cref{eq:rankone-eigenvectors}, note that when $\rho=1/2$, we have $2(1-\rho)^{-1/2}\leqslant 4\leqslant 4dV^2$, and $2\rho^{-1}(d-i+1)V^2\leqslant 4 dV^2$, and when $\maxtwo{\lambda_i}{\lambda_j}/\mintwo{\lambda_i}{\lambda_j}>1+4dV^2$, by \Cref{eq:eigenvectors-rankone-lambdai>lambdaj-distinct},
		 then,
		 $$(1-(1+dV^2)\mintwo{\lambda_i}{\lambda_j}/\maxtwo{\lambda_i}{\lambda_j})^{-1}\leqslant (1+4dV^2)/(4dV^2)\leqslant 5/4,
		 $$
		 as $V\geqslant d^{-1/2}$, and thus \Cref{eq:rankone-eigenvectors} holds. 
		If otherwise $\maxtwo{\lambda_i}{\lambda_j}\leqslant (1+4dV^2)\mintwo{\lambda_i}{\lambda_j}$, as $\maxtwo{\lambda_i}{\lambda_j}/\mintwo{\lambda_i}{\lambda_j} \geqslant 1$, we find $|[e_i]_j|\leqslant 1 \leqslant (1+4dV^2)\sqrt{\maxtwo{\lambda_i}{\lambda_j}/\mintwo{\lambda_i}{\lambda_j}} 
		$. Since $1\leqslant dV^2$, then $(1+4dV^2) \leqslant 5d V^2$ and \Cref{eq:rankone-eigenvectors} holds.
	\end{proof}

	\subsection{Sum of $m$ rank-one matrices pertubation} \label{sec-rankm}
	Our final
	\Cref{th:eigenvectors-rankm} generalizes \Cref{p:eigenvectors-of-rank-one-update-distincts-eigenvalues} to any value $m\geqslant1$ and is obtained by induction using \Cref{eq:eigenvectors-rankone-lambdai>lambdaj-distinct}. \Cref{th:eigenvectors-rankm} implies in particular \Cref{th:eigenvectors-rankm-symmetric} via the spectral theorem.
	\begin{theorem}
		Let $D=\diag(\lambda_1,\dots,\lambda_d)$ be a diagonal matrix with $\lambda_1\geqslant\dots\geqslant\lambda_d>0$. Let $V\geqslant 1/\sqrt d$ and consider a sequence of vectors $v^{(i)}\in\bR^d$ such that $\|v^{(i)}\|_\infty\leqslant V$ for all $i\in\bN$. For $m\in\bN$, let 
		$
			A^{(m)} {}={} D +\sqrt{D} \sum_{i=1}^m [v^{(i)}][v^{(i)}]^T \sqrt{D}
		$
		and $\nu^{(m)}_1\geqslant\dots\geqslant\nu^{(m)}_d$ the eigenvalues of $A^{(m)}$ and $(e_1^{(m)},\dots,e_d^{(m)})$ a corresponding orthonormal system of eigenvectors.
		Then
		\begin{equation}
			|[e_i^{(m)}]_j| \leqslant C_m\sqrt{\frac{\mintwo{\lambda_i}{\lambda_j}}{\maxtwo{\lambda_i}{\lambda_j}}}
			\quad\text{for all~$i,j\in\{1,\dots,d\}$ and $m\in\bN$}
			\label{eq:rank-mu-eigenvectors}
		\end{equation}
		with
		$C_0=1$ and $C_{m+1} =   5d^7V^4C_m^5\sqrt{1+dmV^2}$.
		\label{th:eigenvectors-rankm}
	\end{theorem}
	
	\newcommand{\minmax}{\alpha}
	
	\begin{proof}
		For $a,b>0$, denote $\minmax(a,b)=\sqrt{\mintwo{a}{b}/\maxtwo{a}{b}}$. 
		Let $m\in\bN$ and assume that \Cref{eq:rank-mu-eigenvectors} holds which is true if $m=0$ since $C_0=1$. Observe now that
		$
			A^{(m+1)} = A^{(m)} + \sqrt{D} [v^{(m+1)}][v^{(m+1)}]^T\sqrt{D} .
		$
		In the system of coordinates $\mathcal{B}^{(m)}\coloneqq(e_1^{(m)},\dots,e_d^{(m)})$, $A^{(m)}$ writes as $D^{(m)}\coloneqq\diag(\nu_1^{(m)},\dots,\nu_d^{(m)})$. Since $\lambda_1,\dots,\lambda_d>0$, and as $A^{(m)}\succeq D$, then $\nu_i^{(m)}\geqslant\lambda_i>0$ for $i\in\{1,\dots,d\}$, and 
		$$
			\left\langle\sqrt{D}v^{(m+1)},e_i^{(m)}\right\rangle = \sum_{j=1}^d [e_i^{(m)}]_j \sqrt{\lambda_j}[v^{(m+1)}]_j  
			 = \sqrt{D^{(m)}_{ii}} \times  \sum_{j=1}^d \sqrt{{\lambda_j}/{\nu_i^{(m)}}} [e_i^{(m)}]_j
			 [v^{(m+1)}]_j  
			\eqqcolon \left[\sqrt{D^{(m)}} w^{(m+1)}\right]_i .
		$$
		Hence
		$
			\left[ A^{(m+1)} \right]_{\mathcal{B}^{(m)}} = D^{(m)} + \sqrt{D^{(m)}} [w^{(m+1)}]$ $[w^{(m+1)}]^T \sqrt{D^{(m)}}
		$
		with 
		$$
			\|w^{(m+1)}\|_\infty \leqslant \sum_{j=1}^d\sqrt{\fracp{\lambda_j}{\nu_i^{(m)}}}|[e_i^{(m)}]_j| \times V \leqslant \sum_{j=1}^d\sqrt{\fracp{\lambda_j}{\lambda_i}}|[e_i^{(m)}]_j| \times V \leqslant  dV \times C_m.
		$$
		We apply \Cref{p:eigenvectors-of-rank-one-update-distincts-eigenvalues} 
		to $\left[ A^{(m+1)} \right]_{\mathcal{B}^{(m)}}$ so that, for $i,k\in\{1,\dots,d\}$,
		$$
			|\langle e_i^{(m+1)},e_k^{(m)}\rangle |= |[ [e_i^{(m+1)}]_{\mathcal{B}^{(m)}}]_k| \leqslant 5d^6V^4C_m^4 \minmax(\nu_i^{(m)},\nu_k^{(m)}) .
		$$
		By \Cref{l:eigenvalues-rankm},
		$
			|\langle e_i^{(m+1)},e_k^{(m)}\rangle | \leqslant (5d^6V^4C_m^4)(1+dmV^2)^{1/2}$ $\minmax(\lambda_i,\lambda_k)
		$.
		Since 
		$
		|[e_k^{(m)}]_j|\leqslant C_m\minmax(\lambda_i,\lambda_k)
		$,
		then, 
		$$
			|[e_i^{(m+1)}]_j|  \leqslant \sum_{k=1}^d |[e_k^{(m)}]_j|\times |\langle e_i^{(m+1)},e_k^{(m)}\rangle | 
			 \leqslant d C_m \times 5d^6V^4C_m^4(1+dmV^2)^{1/2}\minmax(\lambda_i,\lambda_j)
			 .
		$$
		This proves by induction that \Cref{eq:rank-mu-eigenvectors} holds for all $m\in\bN$.
	\end{proof}
	
	\subsection{Thightness}
	
	\begin{figure}[t]
		\begin{minipage}[c]{0.6\textwidth}
			\begin{subfigure}[c]{0.4\textwidth}
				\includegraphics[height=0.188\textheight]{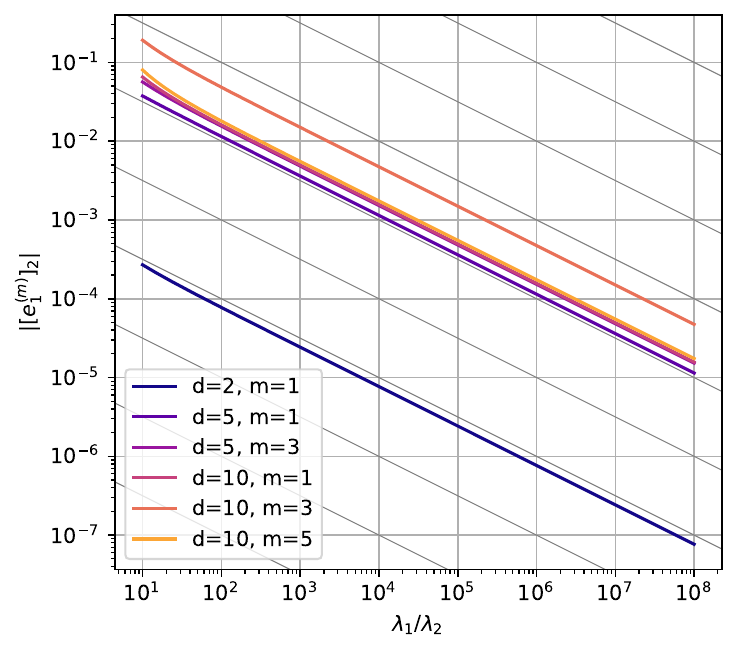}
				\caption{$j=2$}
			\end{subfigure}\hspace{2.0em}
			\begin{subfigure}[c]{0.4\textwidth}
				\includegraphics[height=0.188\textheight]{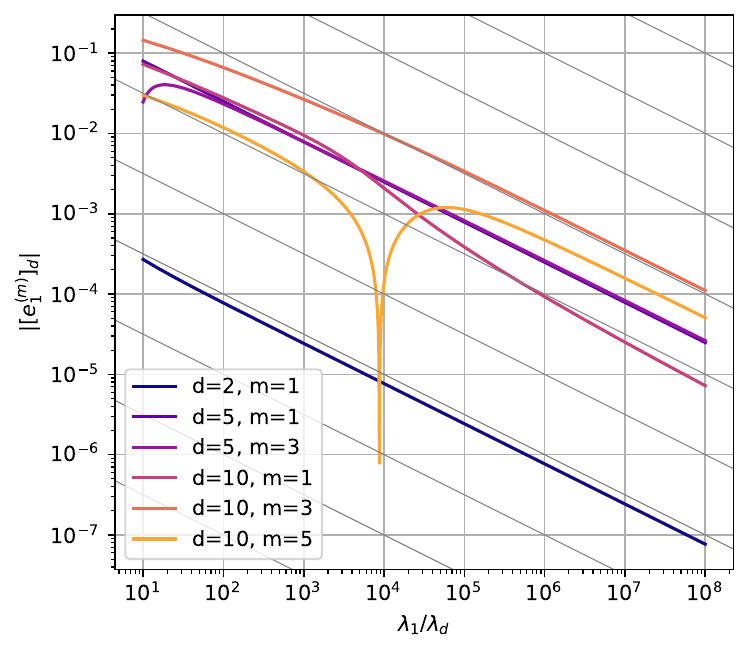}
				\caption{$j=d$}
			\end{subfigure}\hspace{-3.0em}
		\end{minipage}
		\begin{minipage}[c]{0.4\textwidth}
			\caption{Value of $|[e_1^{(m)}]_j|$ as a function of $\lambda_1/\lambda_j$ where $e_1^{(m)}$ is an eigenvector associated to the largest eigenvalue of $A^{(m)}$ 
			from \eqref{eq:diagonal-plus-rankm} for different dimensions and values of $m$ as given in the legend.
The $v^{(i)}$ are independent standard Gaussian vectors (with the same realization for all values of $\lambda_1$) and the eigenvalues of the diagional matrix $B$ are chosen uniformly on a $\log$ scale between $\lambda_d=1$ and $\lambda_1$. The value $|[e_1^{(m)}]_j|$ behaves consistent with $\Theta(\sqrt{\lambda_j/\lambda_1})$. }
			\label{fig:decrease-eigenvectors}
		\end{minipage}
	\end{figure}

	\Cref{fig:decrease-eigenvectors} shows numerical computations of coordinates of the first eigenvector of $A^{(m)}$ in dimension $2$, $5$, $10$.
	The coordinates seem to obey $\Theta(\mintwo{\lambda_i}{\lambda_j}/\maxtwo{\lambda_i}{\lambda_j})$
	in all cases which suggests that this rate in our upper bounds is tight. However we do not expect the constant $C_m$ given in \Cref{th:eigenvectors-rankm} to be tight.

	\subsubsection*{Acknowledgments}
	The authors would like to thank Stéphane Gaubert for his constructive remarks and feedback on a prior version of the manuscript as well as the anonymous referee for their valuable review and comments.
	
	\bibliographystyle{elsarticle-num} 
	\bibliography{biblio-rankm.bib}
	
	
		
		
		
\end{document}